# Square-free Integers and Infinite products

**M. Ramesh Kumar**
Phone No: +91 9840913580
Email address: ramjan_80@yahoo.com
Home page url: http//ramjan07.page.tl/

**Abstract**

We study few properties of square-free integers in certain equations. Using this property, we derive some infinite products in powers of square free numbers. Also, we present a method, to convert power series and trigonometric series to infinite products. Infinite products of few elementary trigonometric functions and factorials for large numbers are shown as examples.

**Introduction**

In Number theory, a square-free integer is, one divisible by no perfect square, except 1. For example, 10 is square-free but 18 is not, as it is divisible by $9 = 3^2$. The smallest square-free integers are 1, 2, 3, 5, 6, 7, 10, 11, 13, 14, 15 etc. Generally, these integers are very less in applications. In this paper, we present few properties of these integers and also we present some infinite products in terms of these integers.

**Definition**

If $n$ is square free Integer, then we can write $n = n_1 n_2 n_3 \ldots n_k$, then $k$ is number of primes in '$n$'. Therefore, $k$ is order of square free Integer $n$.

1. If $|x| < 1$, Let us take

$$x = \sum_{k=1}^{\infty} a_k \log(1 - x^k) \tag{1}$$

So that

$$e^x = \prod_{k=1}^{\infty} \left(1 - x^k\right)^{a_k} \tag{2}$$

$$x = -\sum_{k=1}^{\infty} a_k \left( \frac{x^k}{1} + \frac{x^{2k}}{2} + \frac{x^{3k}}{3} + \ldots \infty \right) \tag{3}$$

$$-x = x\left(\frac{a_1}{1}\right) + x^2\left(\frac{a_1}{2} + a_2\right) + x^3\left(\frac{a_1}{3} + a_3\right) + x^4\left(\frac{a_1}{4} + \frac{a_2}{2} + a_4\right) + \ldots \infty \tag{4}$$

Now, comparing coefficients

$$a_1 = -1$$
$$\frac{a_1}{2} + a_2 = 0$$
$$\frac{a_1}{3} + a_3 = 0 \qquad (5)$$
$$\frac{a_1}{4} + \frac{a_2}{2} + a_4 = 0$$
$$\ldots\ldots\ldots\ldots$$

For example, the $36^{th}$ equation will be in the following form

$$\frac{a_1}{36} + \frac{a_2}{18} + \frac{a_3}{12} + \frac{a_4}{9} + \frac{a_6}{6} + \frac{a_9}{4} + \frac{a_{12}}{3} + \frac{a_{18}}{2} + a_{36} = 0 \qquad (6)$$

30 is square free integer, its factors are 1, 2, 3, 5, 6,10,15,30. Then, the $30^{th}$ equation will be in the form

$$\frac{a_1}{30} + \frac{a_2}{15} + \frac{a_3}{10} + \frac{a_5}{6} + \frac{a_{10}}{3} + \frac{a_{15}}{2} + a_{30} = 0$$

Solving all equations in (5), we find that

$$a_n = \frac{(-1)^{k+1}}{n} \qquad \text{if } n \text{ is square free}$$
$$a_n = 0 \qquad \text{Otherwise} \qquad (7)$$

Where $n_1, n_2, n_3, \ldots, n_k$ are distinct primes

Replacing $x$ by $-x$

$$-x = a_k \sum_{k=1}^{\infty} \log(1 - (-x)^k) \qquad (8)$$

Then,

$$e^{-x} = \prod_{k=1}^{\infty} (1 - (-x)^k)^{a_k} \qquad (9)$$

Subtracting (1) and (8)

$$2x = a_1 \log\left|\frac{1-x}{1+x}\right| + a_3 \log\left|\frac{1-x^3}{1+x^3}\right| + a_5 \log\left|\frac{1-x^5}{1+x^5}\right| + \ldots + \infty \qquad (10)$$

(ie)

$$e^{2x} = \prod_{k=1}^{\infty} \left( \frac{1 - x^{2k-1}}{1 + x^{2k-1}} \right)^{a_{2k-1}} \tag{11}$$

Putting $x = \frac{1}{2}$ in (2)

$$\sqrt{e} = \prod_{k=1}^{\infty} \left(1 - \frac{1}{2^k}\right)^{a_k} \tag{12}$$

Putting $x = \frac{1}{2}$ in (11)

$$e = \prod_{k=1}^{\infty} \left( \frac{2^{2k-1} - 1}{2^{2k-1} + 1} \right)^{a_{2k-1}} \tag{13}$$

2. For $x^2 \leq 1$ and $x\cos\theta \neq 1$, Let us take

$$\sum_{k=1}^{\infty} a_k \log\left(1 - 2x^k \cos k\theta + x^{2k}\right) \tag{14}$$

Where, all $a$'s are defined from (7)

$$= \sum_{k=1}^{\infty} a_k \left( -2 \sum_{i=1}^{\infty} \frac{\cos i\theta}{i} x^{ki} \right)$$

Simplifying, we get

$$= -2\left( x\cos\theta a_1 + x^2 \cos 2\theta \left(\frac{a_1}{2} + a_2\right) + x^3 \cos 3\theta \left(\frac{a_1}{3} + a_3\right) + x^4 \cos 4\theta \left(\frac{a_1}{4} + \frac{a_2}{2} + a_4\right) + \ldots \infty \right)$$

Then, using equations (5),

$$= 2x\cos\theta$$

Thus,

$$2x\cos\theta = \sum_{k=1}^{\infty} a_k \log\left(1 - 2x^k \cos k\theta + x^{2k}\right) \tag{15}$$

Therefore, we get,

$$e^{2x\cos\theta} = \prod_{k=1}^{\infty} \left(1 - 2x^k \cos k\theta + x^{2k}\right)^{a_k} \tag{16}$$

Replacing $x$ by $-x$

$$-2x\cos\theta = \sum_{k=1}^{\infty} a_k \log\left(1 - 2(-x^k)\cos k\theta + x^{2k}\right) \tag{17}$$

Subtracting (17) from (15)

$$4x\cos\theta = \sum_{k=1}^{\infty} a_{2k-1} \log\left|\frac{1 + 2x^{2k-1}\cos\overline{2k-1}\theta + x^{4k-2}}{1 - 2x^{2k-1}\cos\overline{2k-1}\theta + x^{4k-2}}\right|$$

$$e^{4x\cos\theta} = \prod_{k=1}^{\infty}\left(\frac{1 + 2x^{2k-1}\cos\overline{2k-1}\theta + x^{4k-2}}{1 - 2x^{2k-1}\cos\overline{2k-1}\theta + x^{4k-2}}\right)^{a_{2k-1}} \tag{18}$$

If $x=1$, then (16) gives

$$e^{2\cos\theta} = \prod_{k=1}^{\infty}\left(4\sin\frac{k\theta}{2}\right)^{2a_k} \tag{19}$$

Replacing $\theta$ by $\pi - \theta$

$$e^{-2\cos\theta} = \prod_{k=1}^{\infty}\left(4\sin\frac{k(\pi-\theta)}{2}\right)^{2a_k} \tag{20}$$

If $x=1$, then (18) gives

$$e^{4\cos\theta} = \prod_{k=1}^{\infty} \tan^{2a_{2k-1}}(2k-1)\frac{\theta}{2} \tag{21}$$

Adding (15) and (17)

$$\prod_{k=1}^{\infty}\left(\left(1 - 2x^{2k-1}\cos\overline{2k-1}\theta + x^{4k-2}\right)\left(1 + 2x^{2k-1}\cos\overline{2k-1}\theta + x^{4k-2}\right)\right)^{a_{2k-1}}$$
$$= \frac{1}{\prod_{k=1}^{\infty}\left(1 + 2x^{2k}\cos 2k\theta + x^{4k}\right)^{2a_{2k}}} \tag{22}$$

3. Let us take, if $|x| \leq 1$,

$$x = \sum_{k=1}^{\infty} b_k \log(1 + x^k) \tag{23}$$

So that,

$$e^x = \prod_{k=1}^{\infty}\left(1 + x^k\right)^{b_k} \tag{24}$$

$$x = \sum_{k=1}^{\infty} b_k \left( \frac{x^k}{1} - \frac{x^{2k}}{2} + \frac{x^{3k}}{3} - \ldots \infty \right)$$

$$x = x\left(\frac{b_1}{1}\right) + x^2\left(-\frac{b_1}{2} + b_2\right) + x^3\left(\frac{b_1}{3} + b_3\right) + x^4\left(-\frac{b_1}{4} - \frac{b_2}{2} + b_4\right) + \ldots \infty$$

Now, comparing coefficients

$$\left.\begin{array}{l} b_1 = 1 \\ -\dfrac{b_1}{2} + b_2 = 0 \\ \dfrac{b_1}{3} + b_3 = 0 \\ -\dfrac{b_1}{4} - \dfrac{b_2}{2} + b_4 = 0 \\ \ldots\ldots\ldots\ldots\ldots\ldots \end{array}\right\} \tag{25}$$

For example, the 36$^{th}$ equation will be in the following form

$$-\frac{b_1}{36} - \frac{b_2}{18} - \frac{b_3}{12} + \frac{b_4}{9} - \frac{b_6}{6} + \frac{b_9}{4} + \frac{b_{12}}{3} - \frac{b_{18}}{2} + b_{36} = 0 \tag{26}$$

Then, solving above equations we get,

$$\left.\begin{array}{ll} b_1 = 1, b_2 = \frac{1}{2} & \\ b_m = -\dfrac{1}{m} & \text{if } m \text{ is odd prime} \\ b_{n_1 n_2 \ldots n_k} = b_{n_1} b_{n_2} \ldots b_{n_k} & \text{if } n \text{ is square free} \\ b_{2^k m} = \dfrac{1}{2} b_m & \text{if } k \in N, m \text{ is square free} \\ b_n = 0 & \text{Otherwise} \end{array}\right\} \tag{27}$$

Hence, It is easily find that

$$\sum_{k=1}^{\infty} b_k \log\left(1 + 2x^k \cos k\theta + x^{2k}\right) \tag{28}$$

$$= \sum_{k=1}^{\infty} b_k \left( -2 \sum_{i=1}^{\infty} \frac{\cos i\theta}{i} \left(-x^k\right)^i \right) \tag{29}$$

$$= 2\left( x \cos \theta b_1 + x^2 \cos 2\theta \left(-\frac{b_1}{2} + b_2\right) + x^3 \cos 3\theta \left(-\frac{b_1}{3} + b_3\right) + x^4 \cos 4\theta \left(-\frac{b_1}{4} - \frac{b_2}{2} + b_4\right) + \ldots \infty \right)$$

$$= 2x\cos\theta$$

So that,

$$2x\cos\theta = \sum_{k=1}^{\infty} b_k \log\left(1 + 2x^k \cos k\theta + x^{2k}\right) \qquad (30)$$

Therefore, we get,

$$e^{2x\cos\theta} = \prod_{k=1}^{\infty} \left(1 + 2x^k \cos k\theta + x^{2k}\right)^{b_k} \qquad (31)$$

Put $x = 1$ in (24), we find that

$$\sum_{k=1}^{\infty} b_k = \frac{1}{\log 2} \qquad (32)$$

4. Let $p(x)$ be a power series in '$x$', as follow as

$$p(x) = c_1 x + c_2 x^2 + c_3 x^3 + c_4 x^4 + \ldots \infty \qquad (33)$$

Using equation (1), we can write,

$$= c_1 \left(a_1 \log(1-x) + a_2 \log(1-x^2) + a_3 \log(1-x^3) + \ldots \infty\right)$$
$$+ c_2 \left(a_1 \log(1-x^2) + a_2 \log(1-x^4) + a_3 \log(1-x^6) + \ldots \infty\right)$$
$$+ c_3 \left(a_1 \log(1-x^3) + a_2 \log(1-x^6) + a_3 \log(1-x^9) + \ldots \infty\right) + \ldots \infty$$

Rearranging, we get

$$= c_1 a_1 \log(1-x) + (c_1 a_2 + c_2 a_1)\log(1-x^2) + (c_1 a_3 + c_3 a_1)\log(1-x^3) + \ldots \infty$$

Let us take

$$p_k = c_1 a_k + c_2 a_{k/2} + c_3 a_{k/3} + \ldots + c_l a_{k/l} + \ldots + c_k a_1 \qquad (34)$$

Where $1, 2, 3, \ldots, l, \ldots, k$ are the factors of $k$

So that,

$$p(x) = p_1 \log(1-x) + p_2 \log(1-x^2) + p_3 \log(1-x^3) + \ldots \infty$$

$$e^{p(x)} = (1-x)^{p_1}(1-x^2)^{p_2}(1-x^3)^{p_3}\ldots\infty \qquad (35)$$

If we use (15)

$$e^{2p(x)} = \prod_{k=1}^{\infty} (1 - 2x^k \cos k\theta + x^{2k})^{p_1} \tag{36}$$

Where

$$p_k = \frac{c_1}{\cos\theta} a_k + \frac{c_2}{\cos 2\theta} a_{k/2} + \frac{c_3}{\cos 3\theta} a_{k/3} + \ldots + \frac{c_l}{\cos l\theta} a_{k/l} + \ldots + \frac{c_k}{\cos k\theta} a_1 \tag{37}$$

Similarly from (23), we can easily find,

$$e^{p(x)} = (1+x)^{p_1'} (1+x^2)^{p_2'} (1+x^3)^{p_3'} \ldots \infty \tag{38}$$

Where,

$$p_k' = c_1 b_k + c_2 b_{k/2} + c_3 b_{k/3} + \ldots + c_l b_{k/l} + \ldots + c_k b_1 \tag{39}$$

From (30), we find that

$$e^{2p(x)} = \prod_{k=1}^{\infty} (1 + 2x^k \cos k\theta + x^{2k})^{p_1'} \tag{40}$$

Where

$$p_k' = \frac{c_1}{\cos\theta} b_k + \frac{c_2}{\cos 2\theta} b_{k/2} + \frac{c_3}{\cos 3\theta} b_{k/3} + \ldots + \frac{c_l}{\cos l\theta} b_{k/l} + \ldots + \frac{c_k}{\cos k\theta} b_1 \tag{41}$$

Let $q(x)$ be a power series in $x$, in odd degree as follow as

$$q(x) = d_1 x + d_3 x^3 + d_5 x^5 + \ldots \infty \tag{42}$$

$$e^{2q(x)} = \left(\frac{1-x}{1+x}\right)^{q_1} \left(\frac{1-x^3}{1+x^3}\right)^{q_3} \left(\frac{1-x^5}{1+x^5}\right)^{q_5} \ldots \infty \tag{43}$$

Let us take

$$q_k = d_1 a_k + d_3 a_{k/3} + \ldots + d_l a_{k/l} + \ldots + d_k a_1 \tag{44}$$

Where $1, 3, \ldots, l, \ldots, k$ are the odd factors of $k$

$$e^{4q(x)} = \prod_{k=1}^{\infty} \left(\frac{1 - 2x^k \cos k\theta + x^{2k}}{1 + 2x^k \cos k\theta + x^{2k}}\right)^{q_k} \tag{45}$$

$$q_k = \frac{c_1}{\cos\theta} a_k + \frac{c_3}{\cos 3\theta} a_{k/3} + \ldots + \frac{c_l}{\cos l\theta} a_{k/l} + \ldots + \frac{c_k}{\cos k\theta} a_1 \tag{46}$$

5. We can convert Trigonometric series to infinite products

$$p(\theta) = c_1 \cos\theta + c_2 \cos 2\theta + c_3 \cos 3\theta + c_4 \cos 4\theta + \ldots \infty \tag{47}$$

$$e^{2p(\theta)} = \prod_{k=1}^{\infty} \left(1 - 2x^k \cos k\theta + x^{2k}\right)^{p_1} \tag{48}$$

Where

$$p_k = \frac{c_1}{x} a_k + \frac{c_2}{x^2} a_{k/2} + \frac{c_3}{x^3} a_{k/3} + \ldots + \frac{c_l}{x^l} a_{k/l} + \ldots + \frac{c_k}{x^k} a_1 \tag{49}$$

$$e^{2p(\theta)} = \prod_{k=1}^{\infty} (1 + 2x^k \cos k\theta + x^{2k})^{p_1'} \tag{50}$$

Where

$$p_k' = \frac{c_1}{x} b_k + \frac{c_2}{x^2} b_{k/2} + \frac{c_3}{x^3} b_{k/3} + \ldots + \frac{c_l}{x^l} b_{k/l} + \ldots + \frac{c_k}{x^k} b_1 \tag{51}$$

$$q(\theta) = d_1 \cos\theta + d_3 \cos 3\theta + d_5 \cos 5\theta + \ldots \infty \tag{52}$$

$$e^{4q(\theta)} = \prod_{k=1}^{\infty} \left(\frac{1 - 2x^k \cos k\theta + x^{2k}}{1 + 2x^k \cos k\theta + x^{2k}}\right)^{q_k} \tag{53}$$

Where $1, 3, \ldots, l, \ldots, k$ are the odd factors of $k$

$$q_k = \frac{d_1}{x} a_k + \frac{d_3}{x^3} a_{k/3} + \ldots + \frac{d_l}{x^l} a_{k/l} + \ldots + \frac{d_k}{x^k} a_1 \tag{54}$$

Similarly, we can express

$$p(x, \theta) = c_1 x \cos\theta + c_2 x^2 \cos 2\theta + c_3 x^3 \cos 3\theta + c_4 x^4 \cos 4\theta + \ldots \infty \tag{55}$$

6. We can express elementary trigonometry functions in infinite products; few of them are listed below

$$\log\left(\frac{x}{\sin x}\right) = \frac{2B_2}{1.2!} x^2 + \frac{2^3 B_4}{2.4!} x^4 + \frac{2^5 B_6}{3.6!} x^6 + \ldots + \infty, \quad 0 < x < 1$$

$$\frac{x}{\sin x} = (1 - x^2)^{p_1} (1 - x^4)^{p_2} (1 - x^6)^{p_3} \ldots \infty \tag{56}$$

$$p_k = c_1 a_k + c_2 a_{k/2} + c_3 a_{k/3} + \ldots + c_l a_{k/l} + \ldots + c_k a_1, \text{ Where } c_k = \frac{2^{2k-1} B_{2k}}{k.2k!}$$

$$\log\frac{1}{\cos x} = \frac{2(2^2 - 1)|B_2|}{1.2!} x^2 + \frac{2^3(2^4 - 1)|B_4|}{2.4!} x^4 + \frac{2^5(2^6 - 1)|B_6|}{3.6!} x^6 + \ldots + \infty$$

$$\frac{1}{\cos x} = (1-x^2)^{p_1}(1-x^4)^{p_2}(1-x^6)^{p_3}\ldots\infty \tag{57}$$

$$p_k = c_1 a_k + c_2 a_{k/2} + c_3 a_{k/3} + \ldots + c_l a_{k/l} + \ldots + c_k a_1, \text{ Where } c_k = \frac{2^{2k-1}(2^k - 1)|B_{2k}|}{k.2k!}$$

Another interesting example is, factorial of large number, it is known that

$$\log\left(\frac{(n-1)!}{\sqrt{2\pi}n^{n-\frac{1}{2}}e^{-n}}\right) = \frac{B_2}{1.2}\frac{1}{n} + \frac{B_4}{3.4}\frac{1}{n^3} + \frac{B_6}{5.6}\frac{1}{n^5} + \ldots + \infty$$

$$\left(\frac{(n-1)!}{\sqrt{2\pi}n^{n-\frac{1}{2}}e^{-n}}\right)^2 = \left(\frac{n-1}{n+1}\right)^{q_1}\left(\frac{n^3-1}{n^3+1}\right)^{q_3}\left(\frac{n^5-1}{n^5+1}\right)^{q_5}\ldots\infty \tag{58}$$

Provided that,

$$q_k = d_1 a_k + d_3 a_{k/3} + \ldots + d_l a_{k/l} + \ldots + d_k a_1, \text{ Where } d_j = \frac{B_{2j}}{2j(2j-1)}$$

7. Now, if we define $s > 1$

$$\Phi(x) = \frac{x}{1^s} + \frac{x^2}{2^s} + \frac{x^3}{3^s} + \ldots + \infty \tag{59}$$

Let us take

$$x = \sum_{k=1}^{\infty} a_k \Phi(x^k) \tag{60}$$

$$-x = x\left(\frac{a_1}{1^s}\right) + x^2\left(\frac{a_1}{2^s} + a_2\right) + x^3\left(\frac{a_1}{3^s} + a_3\right) + x^4\left(\frac{a_1}{4^s} + \frac{a_2}{2^s} + a_4\right) + \ldots \infty$$

Now, comparing coefficients

$$\left.\begin{array}{l} a_1 = -1 \\ \dfrac{a_1}{2^s} + a_2 = 0 \\ \dfrac{a_1}{3^s} + a_3 = 0 \\ \dfrac{a_1}{4^s} + \dfrac{a_2}{2^s} + a_4 = 0 \\ \ldots\ldots\ldots\ldots \end{array}\right\} \tag{61}$$

For example, the 36th equation will be in the following form

$$\frac{a_1}{36^s} + \frac{a_2}{18^s} + \frac{a_3}{12^s} + \frac{a_4}{9^s} + \frac{a_6}{6^s} + \frac{a_9}{4^s} + \frac{a_{12}}{3^s} + \frac{a_{18}}{2^s} + a_{36} = 0 \tag{62}$$

30 is square free integer, its factors are 1, 2, 3, 5, 6,10,15,30. Then, the 30th equation will be in the form

$$\frac{a_1}{30^s} + \frac{a_2}{15^s} + \frac{a_3}{10^s} + \frac{a_5}{6^s} + \frac{a_{10}}{3^s} + \frac{a_{15}}{2^s} + a_{30} = 0$$

Solving all equations in (62), we find that

$$\left. \begin{array}{l} a_n = \dfrac{(-1)^{k+1}}{n^s} \quad \text{if } n \text{ is square free} \\ a_n = 0 \quad \quad \quad \quad \text{Otherwise} \end{array} \right\} \tag{63}$$

Where $n_1, n_2, n_3, \ldots, n_k$ are distinct primes, Put $x=1$ in (60), we find that

$$\sum_{k=1}^{\infty} a_k = \frac{1}{\zeta(s)} \tag{64}$$

if we define $s > 1$

$$\tilde{\Phi}(x) = \frac{x}{1^s} - \frac{x^2}{2^s} - \frac{x^3}{3^s} + \ldots + \infty \tag{65}$$

Let us take

$$x = \sum_{k=1}^{\infty} a_k \tilde{\Phi}(x^k) \tag{66}$$

$$x = x\left(\frac{b_1}{1^s}\right) + x^2\left(-\frac{b_1}{2^s} + b_2\right) + x^3\left(\frac{b_1}{3^s} + b_3\right) + x^4\left(-\frac{b_1}{4^s} - \frac{b_2}{2^s} + b_4\right) + \ldots \infty$$

Now, comparing coefficients

$$\left. \begin{array}{l} b_1 = 1 \\ -\dfrac{b_1}{2^s} + b_2 = 0 \\ \dfrac{b_1}{3^s} + b_3 = 0 \\ -\dfrac{b_1}{4^s} - \dfrac{b_2}{2^s} + b_4 = 0 \\ \ldots\ldots\ldots\ldots\ldots\ldots \end{array} \right\} \tag{67}$$

For example, the 36th equation will be in the following form

$$-\frac{b_1}{36^s}-\frac{b_2}{18^s}-\frac{b_3}{12^s}+\frac{b_4}{9^s}-\frac{b_6}{6^s}+\frac{b_9}{4^s}+\frac{b_{12}}{3^s}-\frac{b_{18}}{2^s}+b_{36}=0$$

Then, solving above equations we get,

$$\left.\begin{array}{ll} b_1=1, b_2=\tfrac{1}{2} & \\ b_m=-\dfrac{1}{m} & \text{if } m \text{ is odd prime} \\ b_{n_1 n_2 \ldots n_k}=b_{n_1} b_{n_2} \ldots b_{n_k} & \text{if } n \text{ is square free} \\ b_{2^k m}=\dfrac{1}{2} b_m & \text{if } k \in N, m \text{ is square free} \\ b_n=0 & \text{Otherwise} \end{array}\right\} \quad (68)$$

Put $x=1$ in (65), we find that

$$\sum_{k=1}^{\infty} b_k = \frac{1}{\left(1-2^{1-s}\right)\zeta(s)} \tag{69}$$

**List of References**